\documentclass[12pt]{amsart}
\usepackage{mathtools,latexsym,amssymb,amsthm,amsmath,enumerate,epsfig,setspace,bbding,color,graphicx,caption,epstopdf,thmtools,scrextend}
\usepackage[all]{xy}
\usepackage[retainorgcmds]{IEEEtrantools}
\usepackage{hyperref}
\usepackage[margin=10pt,font=small,labelfont=bf]{caption}
\usepackage{tikz}
\usepackage{tensor}
\usepackage{microtype}
\usetikzlibrary{matrix,arrows.meta,decorations.pathmorphing,positioning}

\definecolor{darkred}{rgb}{0.5,0,0}
\definecolor{darkgreen}{rgb}{0,0.5,0}
\definecolor{darkblue}{rgb}{0,0,0.5}
\definecolor{darkorange}{rgb}{0.3,0.6,0.2}
\definecolor{darkyellow}{rgb}{0.75,0.75,0.2}

\numberwithin{equation}{section}
\theoremstyle{plain}

\declaretheorem[name={Theorem},style=plain,numberwithin=section]{thm}

\declaretheorem[name={Definition},style=definition,qed=$\diamondsuit$,sibling=thm]{defi}
\declaretheorem[name={Example},style=definition,qed=$\triangle$,sibling=thm]{ex}

\declaretheorem[name={Remark},style=definition,sibling=thm]{rmk}
\declaretheorem[numbered=no,name={Theorem},style=plain]{thm*}
\declaretheorem[numbered=no,name={Acknowledgements},style=definition]{ack}
\declaretheorem[numbered=no,name={Main Theorem},style=plain]{mainthm}

\newcommand{\R}{\mathbb{R}}
\newcommand{\C}{\mathbb{C}}
\newcommand{\T}{\mathbb{T}}
\newcommand{\Z}{\mathbb{Z}}

\newcommand{\X}{\mathfrak{X}}
\newcommand{\g}{\mathfrak{g}}

\newcommand{\pr}{\mathrm{pr}}
\newcommand{\lin}{\mathrm{lin}}

\newcommand{\Lie}{\mathrm{Lie}}
\newcommand{\G}{\mathcal{G}}
\newcommand{\B}{\mathcal{B}}
\newcommand{\D}{\mathcal{D}}
\newcommand{\F}{\mathcal{F}}
\newcommand{\vol}{\mathrm{vol}}
\renewcommand{\DH}{\mathrm{DH}}
\newcommand{\aff}{\mathrm{aff}}
\newcommand{\red}{\mathrm{red}}

\newcommand{\s}{\mathbf{s}}
\renewcommand{\O}{\mathcal{O}}
\renewcommand{\t}{\mathbf{t}}
\let\oldphi\phi \let\phi\varphi \let\varphi\oldphi

\newcommand\twiddle[1]{{\widetilde{#1}}}
\newcommand{\comment}[1]{}
\newcommand{\sslash}{\mathbin{/\mkern-6mu/}}

\newcounter{prof}

\makeatletter
\@addtoreset{claim}{prof}
\makeatother

\title{Duistermaat-Heckman measures for Hamiltonian groupoid actions}
\date{\today}
\author{Luka Zwaan}
\thanks{This work was partially supported by NSF grant DMS-2003223 and by FWO and FNRS under EOS project G0I2222N}
\begin{document}

\begin{abstract}
Consider a source proper, source connected regular symplectic groupoid acting locally freely and effectively in a Hamiltonian way, and assume that the moment map is proper and has connected fibres. In this case there is an associated \emph{Duistermaat-Heckman measure} on the quotient orbifold. We show that this measure is polynomial with respect to the natural \emph{affine measure}.
\end{abstract}

\maketitle

\section{Introduction}

The study of symmetries has a long history in classical mechanics and its mathematical formalisations. A particularly powerful instance of this is the theory of Hamiltonian group actions, where one can perform a ``double reduction'' using the symmetry (the action) and the conserved quantities (the moment map) \cite{MarsdenWeinstein1974SympRed, Meyer1973SympRed}. This classical notion of Hamiltonian actions of Lie groups on symplectic manifolds has been thoroughly studied and many remarkable results have been obtained, such as singular reduction \cite{SjamaarLerman1991SingRed}, convexity \cite{GuilleminSternberg1982Convexity, Atiyah1982Convexity} and localisation \cite{AtiyahBott84Localisation, BerlineVergne82Localisation}. Here we focus on a result from \cite{dh}, where Hamiltonian actions of a torus with proper moment map are studied. For an effective such action $\T\circlearrowright(X,\omega)\xrightarrow{\mu}\mathfrak{t}^*$, the result relates two measures on the (open) set of regular values $O\subset\mu(X)\subset\mathfrak{t}^*$. One is just the Lebesgue measure on $\mathfrak{t}^*$ and the other, now known as the \emph{Duistermaat-Heckman measure} and written $\mu_{\DH}$, is defined as the pushforward of the Liouville measure on $X$ by the moment map. The result then reads
\begin{equation}\label{eq:classdh}
\mu_{\DH}=\vol_{\red}\cdot\mu_{\mathrm{Leb}}
\end{equation}
where $\vol_{\red}$ assigns to $\xi\in\t^*$ the volume of the reduced space at $\xi$. Moreover, it is shown that $\vol_{\red}$ is a polynomial function.

Several variations of Hamiltonian actions have been introduced, such as quasi-Hamiltonian actions \cite{McDuff1988QuasiHam}, Hamiltonian actions of Poisson-Lie groups \cite{Lu1989PLGMomentMap} and group-valued moment maps \cite{ALekseevMalkinMeinrenken1998GroupValuedMomentMap}. A more general notion of Hamiltonian action that unifies the ones mentioned above can be formulated using \emph{symplectic groupoids}. These are the global objects that integrate Poisson manifolds and they have have played a major role in recent advances in Poisson geometry. An action of a symplectic groupoid on a symplectic manifold is called Hamiltonian if it satisfies a certain multiplicativity condition. In this setting, the moment map is a Poisson map to the base of the symplectic groupoid. In this paper, we show that Duistermaat \& Heckman's result (\ref{eq:classdh}) on polynomial measures generalises to a certain reasonable class of Hamiltonian groupoid actions.

In order to explain our main result, recall that a symplectic groupoid is a Lie groupoid $\G\rightrightarrows M$ carrying a multiplicative symplectic form $\Omega\in\Omega^2(\G)$. There is a unique induced Poisson structure $\pi$ on $M$ making the target map $(\G,\Omega)\to(M,\pi)$ into a Poisson map. An action of $(\G,\Omega)\rightrightarrows(M,\pi)$ on a symplectic manifold $(X,\omega)$ along $\mu:X\to M$ is called \emph{Hamiltonian} if the multiplicativity condition $a^*\omega=\pr_1^*\Omega+\pr_2^*\omega\in\Omega^2(\G\tensor*[_{\s}]{\times}{_{\mu}}X)$ holds. In this paper we consider source proper, source connected, regular symplectic groupoids. This means that the source fibers are connected and compact and that the induced foliation on $M$ is regular. A consequence of this is that that leaf space $B=M/\G$ is an \emph{orbifold} and thus there is a nice theory of measures on the leaf space (see \cite{CrainicMestre2019Measures}): they are in some sense ``transverse measures'' on $M$, with respect to the orbit foliation. The symplectic structure on $\G\rightrightarrows M$ induces an integral affine structure on its leaf space, and this in turn induces a measure on it which we refer to as the \emph{affine measure}. It plays the role of the Lebesgue measure on the dual of the Lie algebra in the classical setting. There is also the \emph{Duistermaat-Heckman measure}, which we define as in the classical case as the pushforward of the Liouville measure on $X$ along the moment map and the quotient map to the leaf space. 

\begin{mainthm}
If the action is locally free and effective and the moment map is proper and has connected fibers, the Duistermaat-Heckman measure is equal to a polynomial function times the affine measure.
\end{mainthm}

Just as in the classical case, we have an exact interpretation of the polynomial function. It is, up to integer factors having to do with the isotropy of $\G$, the function that assigns to an orbit the symplectic volume of the orbit times the symplectic volume of the associated reduced space.

In Section \ref{sec:background} we discuss the necessary background on symplectic groupoids and their actions, integral affine geometry and measures on leaf spaces. In Section \ref{sec:results} we discuss our results, stating the Main Theorem and discussing several remarks and examples in Section \ref{subsec:measresults}. Finally, we give the proof of the Main Theorem in Section \ref{subsec:proof}.

\begin{ack}
I would like to thank my PhD advisor Rui Loja Fernandes for his guidance throughout this project and for his feedback while writing this paper.
\end{ack}

\section{Background}\label{sec:background}

In this section, we review the necessary background and set up our notation.

\subsection{Symplectic groupoids \& Hamiltonian actions}

We recall the basic theory of symplectic groupoids \cite{Karasev1986SymplecticGroupoids, Weinstein1987SymplecticGroupoids, CosteDazordWeinstein1987SymmplecticGroupoids} and their Hamiltonian actions \cite{MikamiWeinstein1988HamiltonianGroupoidActions}. We follow the modern formulation of \cite[Chapter 14]{LecturesPoisson2021}.

\begin{defi}\label{def:sympgpd}
A \emph{symplectic groupoid} is a pair $(\G,\Omega)$ consisting of a Lie groupoid $\G\rightrightarrows M$ and a multiplicative symplectic form $\Omega\in\Omega^2(\G)$.
\end{defi}

\begin{rmk}
Recall that a form $\alpha\in\Omega^{\bullet}(\G)$ is \emph{multiplicative} if 
\[ m^*\alpha=\pr_1^*\alpha+\pr_2^*\alpha\in\Omega^{\bullet}(\G^{(2)}).
\]
Here $\G^{(2)}=\G\tensor*[_{\s}]{\times}{_{\t}}\G$ is the space of composable arrows and $m,\pr_i:\G^{(2)}\to\G$ are the multiplication and projection maps, respectively.
\end{rmk}

Any symplectic groupoid $(\G,\Omega)$ induces a Poisson structure on its base: this is the unique Poisson structure $\pi\in\X^2(M)$ such that the target map $\t:(\G,\Omega)\to(M,\pi)$ is a Poisson map.

\begin{defi}\label{def:hamactiongpd}
Let $(\G,\Omega)$ be a symplectic groupoid and $(X,\omega)$ a symplectic manifold. A groupoid action of $\G$ on $X$ along $\mu:X\to M$ is called \emph{Hamiltonian} if
\[ a^*\omega=\pr_1^*\Omega+\pr_2^*\omega\in\Omega^2(\G\tensor*[_{\s}]{\times}{_{\mu}}X).
\]
Here $a,\pr_2:\G\tensor*[_{\s}]{\times}{_{\mu}}X\to X$ and $\pr_1:\G\tensor*[_{\s}]{\times}{_{\mu}}X\to\G$ are the action and projection maps, respectively.
\end{defi}

In this context we call $\mu$ the \emph{moment map}. It is a Poisson map $\mu:(X,\omega)\to(M,\pi)$. We can display the data of a Hamiltonian action in a diagram as follows:

\begin{center}
\begin{tikzpicture}
\matrix(m)[matrix of math nodes,
row sep=3em, column sep=4em,
text height=1.5ex, text depth=0.25ex]
{
(\G,\Omega) & (X,\omega)  \\
(M,\pi) & \\
};
\begin{scope}[every node/.style={midway,auto,font=\small}]
\draw[->] ([xshift=0.56ex]m-1-1.south) -- ([xshift=0.56ex]m-2-1.north);
\draw[->] ([xshift=-0.56ex]m-1-1.south) -- ([xshift=-0.56ex]m-2-1.north);
\draw[->] (m-1-2) -- node {$\mu$} (m-2-1);
\draw[->] ([xshift=-1.5ex,yshift=-0.56ex]m-1-2.west) .. controls ([xshift=1.5ex,yshift=-1.5em]m-1-1.center) and ([xshift=1.5ex,yshift=1.5em]m-1-1.center) .. ([xshift=-1.5ex,yshift=0.56ex]m-1-2.west);
\end{scope}
\end{tikzpicture}
\end{center}

There is an infinitesimal moment map condition for Hamiltonian actions as well. Let $A=\ker(d\t)|_M$ be the Lie algebroid of $\G$. We have an isomorphism $\sigma_{\Omega}:A\to T^*M$ defined by
\[ \alpha\mapsto -\Omega(\alpha,\cdot)|_M
\]
and denoting the infinitesimal action associated to $\alpha\in\Gamma(A)$ by $\alpha^X\in\X(X)$ we have the formula
\begin{equation}\label{eq:infimoment}
i_{\alpha^X}\omega=\mu^*(\sigma_{\Omega}(\alpha)).
\end{equation}
If $\G\rightrightarrows M$ is source connected, equation (\ref{eq:infimoment}) is in fact equivalent to the multiplicativity condition in Definition \ref{def:hamactiongpd}.

For Hamiltonian groupoid actions, there is also the notion of reduced spaces. For $p\in M$, we denote by $\G_p$ the isotropy group at $p$, and by $\g_p=\Lie(\G_p)$ its Lie algebra. Now $\G_p$ acts on $\mu^{-1}(p)$ and the \emph{reduced space} at $p$ is defined as the quotient $X\sslash_p\G:=\mu^{-1}(p)/\G_p$. If $p$ is a regular value of $\mu$ and the action is proper, $X\sslash_p\G$ inherits an orbifold structure. In fact, it is a symplectic orbifold, the symplectic structure $\omega_{\red,p}$ being induced by presymplectic form $\omega|_{\mu^{-1}(p)}$. Of course, if the action is free then $(X\sslash_p\G,\omega_{\red,p})$ is a symplectic manifold.

\begin{ex}\label{ex:classham}
Let $G$ be a Lie group with Lie algebra $\g$ and let $G\circlearrowright(X,\omega)$ be a Hamiltonian action with equivariant moment map $\mu:X\to\g^*$. The Poisson manifold $(\g^*,\pi_{\lin})$ integrates to the symplectic groupoid $(G\ltimes\g^*,-\omega_{\mathrm{can}})$, the action groupoid of the coadjoint action $G\circlearrowright\g^*$. The symplectic form is the canonical one transported by the isomorphism $T^*G\cong G\times\g^*$ induced by left translations. The action of $G$ induces an action of $(G\ltimes\g^*,-\omega_{\mathrm{can}})$ on $(X,\omega)$ with moment map $\mu$:
\[ (g,\mu(x))\cdot x := gx.
\]
This action is easily checked to be Hamiltonian in the sense of Definition \ref{def:hamactiongpd}.
\end{ex}

\begin{ex} \cite{Lu1989PLGMomentMap, WeinsteinXu1992QYBE}
Let $G$ be a 1-connected complete Poisson-Lie group with dual $G^*$, acting in a Hamiltonian way on a symplectic manifold $(X,\omega)$ with moment map $\mu:X\to G^*$. Recall that this means that the action map and the moment map are Poisson maps and that the moment map condition
\[
i_{a^X}\omega = - \mu^*(a^R)
\]
is satisfied for all $a\in\g$, where $a^R\in\Omega^1(G^*)$ is the right-invariant $1$-form on $G^*$ with value $a$ at the identity. 

The action groupoid $G\ltimes G^*\rightrightarrows G^*$ of the right dressing action $G^*\circlearrowleft G$ turned into a left action by inversion admits a symplectic form turning it into a symplectic groupoid integrating $G^*$. As above, this groupoid acts on $X$ along $\mu$ by
\[ (g,\mu(x))\cdot x := gx
\]
and this action is Hamiltonian.

Note that when $G$ is compact there is an equivariant isomorphism $G^*\cong\g^*$ and the symplectic form $\omega\in\Omega^2(X)$ can be perturbed so as to produce a classical Hamiltonian action \cite{Alekseev1997PLGHamAction}. Hence in this case there is no substantial difference between this example and Example \ref{ex:classham}.
\end{ex}

Let us remark here that the above notions generalise in a straightforward way to the setting of quasi-symplectic groupoids. Recall \cite{BursztynCrainicWeinsteinZhu2004IntDirac, Xu2004MomentumMorita} that a \emph{($\varphi$-twisted) quasi-symplectic groupoid} (also called \emph{presymplectic groupoid}) is a triple $(\G,\Omega,\varphi)$ consisting of a Lie groupoid $\G\rightrightarrows M$, a multiplicative $2$-form $\Omega\in\Omega^2(\G)$ and a closed $3$-form $\varphi\in\Omega^3_{\mathrm{cl}}(M)$ such that $d\Omega=\s^*\varphi-\t^*\varphi$ and such that
\[ \ker(\Omega_p)\cap\ker(d_p\s)\cap\ker(d_p\t) = 0
\]
for all $p\in M$. The corresponding structures on the base are $\varphi$-twisted Dirac structures. As introduced in \cite{Xu2004MomentumMorita}, a Hamiltonian action in this context consists of an action of $\G$ along $\mu:X\to M$ and a $2$-form $\omega\in\Omega^2(X)$ such that $d\omega=\mu^*\varphi$ and
\[
a^*\omega=\pr_1^*\Omega+\pr_2^*\omega\in\Omega^2(\G\tensor*[_{\s}]{\times}{_{\mu}}X)
\]
and such that the nondegeneracy condition
\[ \ker(\omega_x) = \{\alpha^X\mid \alpha\in A_{\mu(x)} \cap \ker(\Omega_{\mu(x)}) \}
\]
holds for all $x\in X$. 

It is a future goal to generalise our results to this setting.

\begin{ex} \cite{ALekseevMalkinMeinrenken1998GroupValuedMomentMap, BursztynCrainicWeinsteinZhu2004IntDirac, Xu2004MomentumMorita}
Let $G$ be a Lie group and assume that $\g$ admits an invariant nondegenerate symmetric bilinear form $(\cdot,\cdot)$. Denote by $\theta,\bar{\theta}\in\Omega^1(G,\g)$ the left- and right-invariant Maurer-Cartan forms respectively and write $\varphi = \frac1{12}(\theta,[\theta,\theta]) = \frac1{12}(\bar{\theta},[\bar{\theta},\bar{\theta}]) \in\Omega^3(G)$ for the Cartan $3$-form. Recall that in this context an action on a manifold $X$ equipped with a $2$-form $\omega\in\Omega^2(X)$ is called \emph{quasi-Hamiltonian} with moment map $\mu:X\to G$ if $\omega$ is invariant, $\mu$ is equivariant (with respect to the conjugation action) and the following conditions hold:
\begin{enumerate}
\item $d\omega=-\mu^*\varphi$,
\item $i_{a^X}\omega=\frac12\mu^*(\theta+\bar{\theta},a)$ for all $a\in\g$, 
\item $\ker(\omega_x)=\{a^X(x)\mid a\in\ker(\mathrm{Ad}_{\mu(x)}+1)\}$ for all $x\in X$.
\end{enumerate}
The action groupoid $G\ltimes G\rightrightarrows G$ of the conjugation action carries a $2$-form which turns it into a $(-\varphi)$-twisted quasi-symplectic groupoid, which integrates the Cartan-Dirac structure associated to $-\varphi$. Given an action as above this groupoid acts again by the formula
\[ (g,\mu(x))\cdot x = gx.
\]
Note that when $G=\T$ is abelian, $G\ltimes G\rightrightarrows G$ is actually an untwisted symplectic groupoid: it is the trivial $\T$-bundle over $\T$.
\end{ex}

\subsection{Integral affine structure}\label{subsec:ias}

We recall how the leaf space of a regular, source connected, proper symplectic groupoid inherits an integral affine structure. For details, see \cite[Section 3]{pmct2}.

\begin{defi}
A \emph{transverse integral affine structure} on a regular foliation $(M,\F)$ is a lattice $\Lambda\subset\nu^*(\F)$ in the conormal bundle to $\F$ which is locally spanned by basic closed $1$-forms.
\end{defi}

\begin{rmk}
Equivalently, one can define a transverse integral affine structure as a foliation atlas whose transition functions are integral affine maps.
\end{rmk}

Consider now a regular, source connected, proper symplectic groupoid $(\G,\Omega)\rightrightarrows(M,\pi)$ and denote the associated symplectic foliation by $\F_{\pi}$. We obtain a lattice $\Lambda\subset\nu^*(\F_{\pi})$ as follows:
\begin{enumerate}
\item for each $p\in M$, the kernel of the exponential map $\g_p\to\G_p$ gives a lattice in $\g_p$ and
\item the isomorphism $\sigma_{\Omega}:\g_p\cong\nu_p^*(\F_{\pi})$ allows us to transport it to the conormal space.
\end{enumerate}
One can show that $\Lambda$ is a transverse integral affine structure on the symplectic foliation. We can think of this as an integral affine structure on the leaf space.

\subsection{Measures on leaf spaces}

A theory of measures on differentiable stacks is laid out in \cite{CrainicMestre2019Measures}. We present a selection of the results in the specific case of measures on leaf spaces of source proper, regular groupoids, for which the theory simplifies significantly.

\subsubsection{Orbifolds}

Let us first settle notation and convention on orbifolds.

\begin{defi}\label{def:orbifold}
An \emph{orbifold atlas} on a topological space $B$ is a proper foliation groupoid $\B\rightrightarrows X$ together with a homeomorphism $q:X/\B\to B$. An \emph{orbifold} is a triple $(B,\B\rightrightarrows X,p)$ of such data. An \emph{equivalence} between orbifolds $(B,\B\rightrightarrows X,q)$ and $(B',\B'\rightrightarrows X',q')$ is a Morita equivalence $(\B\rightrightarrows X)\cong (\B'\rightrightarrows X')$ whose induced homeomorphism $B\cong B'$ intertwines $q$ and $q'$.
\end{defi}

\begin{rmk}
Recall that a Lie groupoid $\B\rightrightarrows X$ is \emph{proper} if $(\s,\t):\B\to X\times X$ is a proper map and that it is called a \emph{foliation groupoid} if the isotropy groups $\B_x$ are discrete for all $x\in X$. A proper foliation groupoid is also called an \emph{orbifold groupoid}.
\end{rmk}

In this paper we are aiming to analyse the leaf spaces of (source) proper, regular groupoids. For such $\G\rightrightarrows X$, collecting the connected components of the identity of the isotropy groups yields a bundle of tori $\mathcal{T}(\G)$ and the quotient $\B(\G)=\G/\mathcal{T}(\G)$ is an orbifold groupoid over $X$. In this way, the leaf spaces we are interested in inherit the structure of an orbifold. In other words, the orbifolds we are dealing with are actually quotients of groupoids on the nose. Thus Definition \ref{def:orbifold} simplifies significantly and we can essentially just think of orbifolds \& equivalences as orbifold groupoids (or even proper, regular groupoids) \& Morita equivalences.

\subsubsection{Measures on manifolds}

We now turn to measures. Let us first make precise what we mean by ``measures'' in the context of smooth manifolds.

\begin{defi}
Let $X$ be a smooth manifold. A \emph{measure} $\mu$ on $X$ is a linear functional $\mu:C_c^{\infty}(X)\to\R$ on the space of compactly supported smooth functions, satisfying $\mu(f)\geq0$ for all $0\leq f\in C_c^{\infty}(X)$.
\end{defi}

\begin{rmk}
This is essentially the definition of a \emph{Radon measure}, which makes sense for any locally compact, Hausdorff space (replacing $C_c^{\infty}(X)$ by $C_c(X)$). For manifolds, the two definitions are equivalent.
\end{rmk}

A special class of measures consists of those arising from densities. Measures of this type are often called \emph{geometric measures}. Recall that a \emph{density} on an $n$-dimensional vector space $V$ is a function $\rho$ that assigns to a basis $(v_1,\ldots,v_n)$ a real number, such that
\[ \rho(Av_1,\ldots,Av_n) = |\det(A)|\cdot\rho(v_1,\ldots,v_n)
\]
for $A\in\mathrm{GL}(V)$. A density on a manifold $X$ is a section of the density bundle $\D_{TX}$, which at each point $x\in M$ consists of the space of densities on $T_xX$. We denote the set of densities by $\D(X)$. Any differential form $\alpha\in\Omega^{\mathrm{top}}(X)$ induces a density $|\alpha|\in\D(X)$. There is a canonical integration map
\[ \int_X:\D(X)\to\R
\]
and this integration generalises that of differential forms. A density is called (strictly) positive if it takes (strictly) positive values at all points. Any positive density $\rho\in\D(X)$ induces a measure $\mu_{\rho}$, defined by
\[\mu_{\rho}(f):=\int_Xf\cdot\rho, \ \ \ \ \ f\in C_c^{\infty}(M).
\]

Geometric measures can be pushed forward along proper submersions using fiber integration. For $q:X\to Y$ a proper submersion and $\rho\in\D(X)$, the \emph{pushforward} of $\rho$ along $q$ is denoted $q_!(\rho)\in\D(Y)$. General measures, even between locally compact, Hausdorff spaces, can be pushed forward along any proper map and we use similar notation. For geometric measures, these two notions of pushforward are compatible.

\subsubsection{Measures on leaf spaces}\label{subsubsec:measleaf}

We now fix a source proper, regular groupoid $\G\rightrightarrows X$. We write $B$ for its leaf space, $q:X\to B$ for the quotient map and $\B=\B(\G)$ for the underlying orbifold groupoid. Since $B$ is locally compact and Hausdorff, we have the notion of Radon measure on it. However, it turns out that, similar to the case of manifolds, it is equivalent to define a measure as a positive linear functional on
\[ C_c^{\infty}(B):=\{f\in C_c(B)\mid f\circ q \in C^{\infty}(X)\}.
\]
We will describe two ways of obtaining measures on $B$. The first is simple: since $q$ is a proper map, we can push measures on $X$ forward to $B$. The second way is more complicated: it describes measures induced by \emph{transverse densities}. 

Let us write $A=\mathrm{Lie}(\G)$ for the Lie algebroid of $\G$, $\g\subset A$ for the isotropy subbundle, $T\F\subset TX$ for the foliation induced by $\G$ and $\nu(\F)$ for the associated normal bundle. Note that $\g$ and $\nu(\F)$ are $\G$-representations.

\begin{defi}\label{def:transdens}
A \emph{transverse density} on $\G$ is a $\G$-invariant positive density $\rho_{\nu}\in\Gamma(\D_{\nu(\F)})$.
\end{defi}

To see how this induces a measure on $B$, we need to choose a strictly positive density $\rho_{\F}\in\Gamma(\D_{T\F})$. We then set $\rho_X:=\rho_{\F}\otimes\rho_{\nu}\in\D(X)$ and we can define a measure $\mu_{\rho_{\nu}}$ on $B$ by the formula
\begin{equation}
\int_B f(b)\,d\mu_{\rho_{\nu}}(b):=\int_X\dfrac{f(q(x))}{\iota(x)\cdot\mathrm{vol}(\O_x,\rho_{\F})}\,d\mu_{\rho_X}(x), \ \ \ \ \ f\in C_c^{\infty}(B).
\end{equation}
Here $\O_x$ denotes the orbit of $\G$ through $x\in X$, $\mathrm{vol}(\O_x,\rho_{\F})$ is its volume with respect to the density $\rho_{\F}$ (restricted to the orbit) and $\iota(x)$ is the number of connected components of the isotropy group $\G_x$. This definition does not depend on the choice of $\rho_{\F}$. We also have the ``fiber integration formula''
\begin{equation}\label{eq:measurecomp}
\int_X f(x)\,d\mu_{\rho_X}(x):=\int_B\iota(b)\cdot\left(\int_{\O_b}f(x)\,d\mu_{\rho_{\F}}(x)\right)\,d\mu_{\rho_{\nu}}(x), \ \ \ \ \ f\in C_c^{\infty}(X).
\end{equation} 
Here $\O_b=q^{-1}(b)$ is the orbit associated to $b\in B$ and $\iota:B\to\Z_{\geq1}$ is defined as above.

\begin{rmk}
As mentioned before, these definitions only work in the specific case of source proper, regular groupoids. For the more general theory, and for a detailed account on how it reduces to the above in our case, see \cite{CrainicMestre2019Measures}.
\end{rmk}

\section{Results}\label{sec:results}

From now on, let \begin{center}
\begin{tikzpicture}
\matrix(m)[matrix of math nodes,
row sep=3em, column sep=4em,
text height=1.5ex, text depth=0.25ex]
{
(\G,\Omega) & (X,\omega)  \\
(M,\pi) & \\
B & \\
};
\begin{scope}[every node/.style={midway,auto,font=\small}]
\draw[->] ([xshift=0.56ex]m-1-1.south) -- ([xshift=0.56ex]m-2-1.north);
\draw[->] ([xshift=-0.56ex]m-1-1.south) -- ([xshift=-0.56ex]m-2-1.north);
\draw[->] (m-1-2) -- node {$\mu$} (m-2-1);
\draw[->] ([xshift=-1.5ex,yshift=-0.56ex]m-1-2.west) .. controls ([xshift=1.5ex,yshift=-1.5em]m-1-1.center) and ([xshift=1.5ex,yshift=1.5em]m-1-1.center) .. ([xshift=-1.5ex,yshift=0.56ex]m-1-2.west);
\draw[->] (m-2-1) -- node {$q$} (m-3-1);
\end{scope}
\end{tikzpicture}
\end{center}
be a locally free, effective Hamiltonian action of a source connected, source proper, regular symplectic groupoid and assume that $\mu$ is proper and has connected fibres. We have written $q:M\to B$ for the projection to the leaf space and denote the symplectic foliation by $\F_{\pi}$. With these assumptions we can and will assume without loss of generality that the moment map is surjective.

\subsection{Linear variation}\label{subsec:linvar}

The linear variation theorem in this setting is in some sense not so interesting: essentially, the action locally looks like a classical Hamiltonian torus action and thus the linear variation is given just as in that setting. However, it is still fruitful to analyse the situation in some detail.

Recall that we are studying the variation in cohomology of the symplectic forms $\omega_{\red,p}$ as $p$ varies through $M$. The first thing to note is that there is no variation in leafwise directions: an arrow $g\in\G$ with $\s(g)=p_1$ and $\t(g)=p_2$ induces an isomorphism $\mu^{-1}(p_1)\cong\mu^{-1}(p_2)$ that intertwines the actions of $\G_{p_1}$ and $\G_{p_2}$ and the forms $\omega|_{\mu^{-1}(p_1)}$ and $\omega|_{\mu^{-1}(p_2)}$. In other words, it identifies the reduced spaces at $p_1$ and $p_2$. This shows that the linear variation actually naturally takes place on the leaf space $B$ and that we should study the variation in directions transverse to the symplectic foliation.

To this end, let us fix a transverse integral affine chart $(U,\phi)$ and let us write $T:=\phi^{-1}(\{0\}\times\R^q)$ for the transversal associated to $(U,\phi)$. The chart trivialises the conormal bundle and the lattice $\Lambda\subset\nu^*(\F_{\pi})$ over $U$ and since $\sigma_{\Omega}$ by definition identifies $\G_p^0\cong\nu_p^*(\F_{\pi})/\Lambda_p$ we obtain a trivialisation of $\mathcal{T}(\G)|_U$ as well. In particular this holds over $T$ as well, and thus we can consider the bundle $\mathcal{T}(\G)|_T$ as a single torus $\T$ acting on $\mu^{-1}(T)$. The latter is a symplectic submanifold of $(X,\omega)$ and the map $\twiddle{\mu}:=\phi\circ\mu:\mu^{-1}(T)\to\R^q$ is easily verified to make the resulting data $\T\circlearrowright(\mu^{-1}(T),\omega|_{\mu^{-1}(T)})\xrightarrow{\twiddle{\mu}}\R^q$ into a classical Hamiltonian torus action. Since the action is locally free, we can apply to it the classical linear variation theorem \cite[Theorem 1.1]{dh}. Note that the reduced spaces of this action are finite covers of the original reduced spaces, since $\G_p^0$ acts instead of $\G_p$.

\subsection{The Duistermaat-Heckman measure}\label{subsec:measresults}

We first introduce the two relevant measures, which now live on the leaf space $B$. On the one hand there is the \emph{affine measure} associated to the transverse integral affine structure, which plays the role of the Lebesgue measure in the classical case. On the other hand there is the \emph{Duistermaat-Heckman measure} associated to the Hamiltonian action. We give an explicit formula for the volumes of the reduced spaces of the action, and finally state our Main Theorem which relates the two measures.

\subsubsection{The affine measure}\label{subsec:affmeas}

From Section \ref{subsec:ias} we get the lattice $\Lambda\subset\nu^*(\F_{\pi})$, which induces a transverse density (Definition \ref{def:transdens}) as follows. Pick any local frame $\{\lambda_1,\ldots,\lambda_q\}$ of $\Lambda$, and define (locally)
\[ \rho_{\nu}:=|\lambda_1\wedge\cdots\wedge\lambda_q|.
\]
It is easily checked that this gives a well-defined transverse density $\rho_{\nu}\in\Gamma(\D_{\nu(\F_{\pi})})$. As described in Section \ref{subsubsec:measleaf} we obtain an induced measure on $B$, which we call the \emph{affine measure} and denote by $\mu_{\mathrm{aff}}$.

In this context, there is a nice candidate for the foliated density $\rho_{\F_{\pi}}$. Since the leaves of the foliation are equipped with symplectic forms, we have a ``foliated Liouville form''
\[ \rho_{\F_{\pi}}:=\left|\frac{\omega_{\F_{\pi}}^{\mathrm{top}}}{\mathrm{top}!}\right|.
\]
The density on $M$ that we work with in computations is then
\[ \rho_M:=\left|\frac{\omega_{\F_{\pi}}^{\mathrm{top}}}{\mathrm{top}!}\right|\otimes\rho_{\nu}.
\]

\subsubsection{The Duistermaat-Heckman measure}\label{subsec:dhmeas}

\begin{defi}
The \emph{Duistermaat-Heckman measure} is defined as the pushforward of the Liouville measure on $X$:
\[ \mu_{\mathrm{DH}}:=(q\circ\mu)_*\left(\frac{\omega^{\mathrm{top}}}{\mathrm{top}!}\right).\qedhere
\]
\end{defi}

We can describe this measure in an alternate way, which gives us an opportunity to analyse the quotient $X/\G$ in some more detail as well. Since the action is locally free, $X/\G$ is an orbifold with atlas $\G\ltimes X\rightrightarrows X$. In fact, one can view it as a ``Poisson orbifold'', the Poisson structure appearing on the level of $X$ as the (regular) Dirac structure $\mathbb{L}_{\omega}$ given by
\[ (\mathbb{L}_{\omega})_x=\{v+w+i_w\omega\mid v\in T_x(\G\cdot x),w\in\ker(d_x\mu)\}.
\]
Here we denote by $\G\cdot x$ the orbit of the $\G$-action through $x$.
The leaf space of this Dirac structure can be thought of as the leaf space of $X/\G$. We could also define the Duistermaat-Heckman measure as the pushforward of the Liouville measure to this leaf space. These two definitions are compatible, as we now explain.

The pullback groupoid $\mu^*\G=X\tensor*[_{\mu}]{\times}{_{\t}}\G\tensor*[_{\s}]{\times}{_{\mu}}X$ equipped with the $2$-form $\omega\oplus-\Omega\oplus-\omega$ is a quasi-symplectic groupoid integrating $\mathbb{L}_{\omega}$ and the standard Morita equivalence $\mu^*\G\cong\G$ establishes an isomorphism between the leaf space of $X/\G$ and $B$. Clearly, this isomorphism interwines the two definitions of the Duistermaat-Heckman measure.

\begin{rmk}
The quasi-symplectic integration of $\mathbb{L}_{\omega}$ also induces a transverse integral affine structure and thus we get an ``affine measure'' on the leaf space of $X/\G$ in the same way as before (see Section \ref{subsec:affmeas}). Since the isomorphism induced by $\mu^*\G\cong\G$ preserves the integral affine structure on the leaf spaces, it also intertwines the two affine measures. All in all, there are no issues working solely on the leaf space $B$.
\end{rmk}

\subsubsection{The Main Theorem}

Before we state the theorem, we need one more bit of notation. Recall that we have a function $\iota:B\to\Z_{\geq1}$ associating to some $b\in B$ the number of connected components of $\G_p$ for any $p\in\O_b$. We denote by $\vol:B\to\R$ the function
\[ b\mapsto\iota(b)\cdot\vol(\O_b,\omega_{\F_{\pi}}).
\]
Similarly, we denote by $\vol_{\red}:B\to\R$ the function
\[ b\mapsto\iota(b)\cdot\vol(X\sslash_p\G,\omega_{\red,p}).
\]
Here we take any $p\in\O_b$ and $\vol(X\sslash_p\G,\omega_{\red,p})$ is the volume of the reduced space. Using the theory from Section \ref{subsubsec:measleaf} we obtain an explicit formula for the volume as follows. Taking the top (nonzero) power of $\omega|_{\mu^{-1}(p)}$ yields a transverse measure (see Definition \ref{def:transdens}) and combining this with the Haar measure on $\G_p$, transported by the action to $X$ as a foliated density, gives the formula
\begin{equation}\label{eq:redspacevol}
\vol(X\sslash_p\G,\omega_{\red,p})=\int_{\mu^{-1}(p)} \frac{\rho_{\mathrm{Haar}}}{\iota(p)} \otimes \left|\frac{\big(\omega|_{\mu^{-1}(p)}\big)^{\mathrm{top}}}{\mathrm{top}!}\right|.
\end{equation}
Here $\rho_{\mathrm{Haar}}$ is normalised according to $\G_p^0$, the identity component of the isotropy group, and $\iota(p)$ is the number  of connected components of $\G_p$.

\begin{rmk}
Strictly speaking, the integrand in equation (\ref{eq:redspacevol}) needs to be modified by a function similar to $\iota$, associating to $x\in \mu^{-1}(p)$ the number of connected components of $\G_x$. However, this function is only not equal to $1$ on a set of measure zero so the equation still holds in the above form.
\end{rmk}

\begin{mainthm}
The Duistermaat-Heckman measure is related to the affine measure by the formula
\begin{equation}\label{eq:mainthm}
\mu_{\DH}=\vol\cdot\vol_{\red}\cdot\mu_{\aff}.
\end{equation}
Moreover, $\vol$ and $\vol_{\red}$ are \emph{polynomial} functions on the leaf space.
\end{mainthm}

We prove the theorem in Section \ref{subsec:proof}.

\begin{rmk}
The notion of a polynomial function makes sense on integral affine manifolds (and orbifolds): since integral affine maps preserve polynomials, one can simply require the function to be polynomial in every integral affine (foliation) chart.
\end{rmk}

\begin{rmk}
The appearance of the two functions $\vol$ and $\vol_{\red}$ corresponds to the two times we push forward the Liouville measure: first along $\mu$, which gives the factor $\vol_{\red}$, and then along $q$, which gives the factor $\vol$.
\end{rmk}

\begin{rmk}
Suppose that we drop the assumption that the action is locally free. It follows from Section \ref{subsec:linvar} that the action of $\mathcal{T}(\G)$ is locally just a torus action. Since the action is effective, this implies that the action of $\mathcal{T}(\G)$ is actually \emph{free} on an invariant open dense set $O\subset X$. In particular, the action of $\G$ is locally free on $O$ and since $\mu(O)$ is consequently open and invariant we can apply our main Theorem to the restricted action of $(\G|_{\mu(O)},\Omega)\rightrightarrows (\mu(O),\pi)$ on $(O,\omega|_O)$.
\end{rmk}

\begin{rmk}[Behaviour under Morita equivalence]\label{rmk:morita}
Suppose we have a symplectic Morita equivalence
\begin{center}
\begin{tikzpicture}
\matrix(m)[matrix of math nodes,
row sep=3em, column sep=4em,
text height=1.5ex, text depth=0.25ex]
{
(\G_1,\Omega_1) & (P,\omega_P) & (\G_2,\Omega_2) \\
(M_1,\pi_1) & & (M_2,\pi_2)\\
 & B & \\
};
\begin{scope}[every node/.style={midway,auto,font=\small}]
\draw[->] ([xshift=0.56ex]m-1-1.south) -- ([xshift=0.56ex]m-2-1.north);
\draw[->] ([xshift=-0.56ex]m-1-1.south) -- ([xshift=-0.56ex]m-2-1.north);
\draw[->] ([xshift=0.56ex]m-1-3.south) -- ([xshift=0.56ex]m-2-3.north);
\draw[->] ([xshift=-0.56ex]m-1-3.south) -- ([xshift=-0.56ex]m-2-3.north);
\draw[->] (m-1-2) -- node {$\mu_1$} (m-2-1);
\draw[->] (m-1-2) -- node[swap] {$\mu_2$} (m-2-3);
\draw[->] ([xshift=-1.5ex,yshift=-0.56ex]m-1-2.west) .. controls ([xshift=1.5ex,yshift=-1.5em]m-1-1.center) and ([xshift=1.5ex,yshift=1.5em]m-1-1.center) .. ([xshift=-1.5ex,yshift=0.56ex]m-1-2.west);
\draw[->] ([xshift=1.5ex,yshift=-0.56ex]m-1-2.east) .. controls ([xshift=-1.5ex,yshift=-1.5em]m-1-3.center) and ([xshift=-1.5ex,yshift=1.5em]m-1-3.center) .. ([xshift=1.5ex,yshift=0.56ex]m-1-2.east);
\draw[->] (m-2-1) -- node {$q_1$} (m-3-2);
\draw[->] (m-2-3) -- node[swap] {$q_2$} (m-3-2);
\end{scope}
\end{tikzpicture}
\end{center} 
between regular, source connected, source proper symplectic groupoids. This induces an equivalence between the categories of Hamiltonian spaces of $(\G_1,\Omega_1)\rightrightarrows(M_1,\pi_1)$ and $(\G_2,\Omega_2)\rightrightarrows(M_2,\pi_2)$: given a Hamiltonian action of the former along $(X,\omega)\xrightarrow{\mu}(M,\pi)$, $\G_1$ acts diagonally on the fibre product $X\tensor*[_{\mu}]{\times}{_{\mu_1}}P$ in a free and proper way, and $\G_2$ acts on the quotient $X\times_{\G_1}P$ along the map $\twiddle{\mu}$ induced by $\mu_2\circ\mathrm{pr}_2$ by 
\[ g\cdot[x,y]:=[x,y\cdot g^{-1}].
\]
Moreover, the form $\omega\oplus(-\omega_P)$ descends to a symplectic form on $X\times_{\G_1}P$ and the $\G_2$-action is Hamiltonian. There is clearly an analogous inverse construction starting from a Hamiltonian space of $(\G_2,\Omega_2)\rightrightarrows(M_2,\pi_2)$.

Now, from the above construction it is clear that when we start with a Hamiltonian action of $(\G_1,\Omega_1)$ with the properties described at the start of Section \ref{sec:results}, the corresponding action of $(\G_2,\Omega_2)$ has the same properties. From the measures and functions in the Main Theorem (\ref{eq:mainthm}), we will see that $\mu_{\aff}$ and $\vol_{\red}$ are preserved under Morita equivalence, while $\vol$ and thus also $\mu_{\DH}$ are not.

Let us write $\mu_{\DH}^i,\mu_{\aff}^i,\vol^i,\vol_{\red}^i$, $i=1,2$, for the relevant objects on $B$ induced by the action of $(\G_i,\Omega_i)$. The standard Morita correspondence shows immediately that $\mu_{\aff}^1=\mu_{\aff}^2$, so there is only one relevant affine measure $\mu_{\aff}$. Furthermore, it is easy to see that for any $y\in P$ the reduced spaces over $\mu_1(y)$ and $\mu_2(y)$ are isomorphic, so that also $\vol_{\red}^1=\vol_{\red}^2$. Of course $\vol^1$ and $\vol^2$ are not equal, which is expected since they give \emph{leafwise} information. We see that the difference between $\mu_{\DH}^1$ and $\mu_{\DH}^2$ is given purely by $\vol^1$ and $\vol^2$, which is data dependent purely on the groupoids, not their actions. Putting the above in different words, the Morita equivalence identifies the measures $\mu_*\left(\frac{\omega^{\mathrm{top}}}{\mathrm{top}!}\right)$ and $\twiddle{\mu}_*\left(\frac{\omega^{\mathrm{top}}}{\mathrm{top}!}\right)$ on $M_1$ and $M_2$ respectively. The difference comes only when pushing forward these measures to $B$.

\end{rmk}

\begin{ex}[The classical case]\label{ex:classcase}
Consider a locally free Hamiltonian torus action $\T\circlearrowright(X,\omega)$ with proper moment map $\mu:X\to\mathfrak{t}^*$. The groupoid is now a bundle of tori $\T\ltimes\mathfrak{t}^*\rightrightarrows \mathfrak{t}^*$ and thus $B=\mathfrak{t}^*$ is smooth. The affine measure is Lebesgue measure on $\mathfrak{t}^*\cong\R^q$ and the Duistermaat-Heckman measure is the classical one as in \cite{dh}. The function $\vol$ has constant value $1$ and $\vol_{\red}$ gives the volume of the reduced spaces. Thus our Main Theorem reduces to the classical Duistermaat-Heckman theorem (\ref{eq:classdh}) in this case.
\end{ex}

\begin{ex}
Consider a symplectic toric manifold $\T^n\circlearrowright(F_0,\omega_0)\xrightarrow{\mu_0}\R^n$. It is well known that, with our conventions, we have $\mu_{\DH}=\mu_{\aff}$. The same holds true for the toric actions of symplectic torus bundles of \cite{FernandesMol2024KahlerMetric}. Note that in this case the leaf space is the base of the torus bundle, and thus smooth.

More generally we have $\mu_{\DH}=\vol\cdot\mu_{\aff}$ for the faithful multiplicity-free spaces of \cite{Mol2023MultFree}.
\end{ex}

\begin{ex}[The free case]\label{ex:freecase}
If the action is free, the situation simplifies significantly. In this case, the quotient $X_{\red}:=X/\G$ is a smooth manifold endowed with a Poisson structure $\pi_{\red}$ induced from $\omega$. In fact, $(X_{\red},\pi_{\red})$ is again a Poisson manifold of source proper type: the gauge groupoid $\big((X\tensor*[_{\mu}]{\times}{_{\mu}}X)/\G,\omega\oplus-\omega\big)$ provides a source connected, source proper symplectic integration. Moreover, $(X,\omega)$ gives a symplectic Morita equivalence with $(\G,\Omega)\rightrightarrows(M,\pi)$. It is not hard to show that for any symplectic Morita equivalence between regular, source connected, source proper symplectic groupoids $(\G_1,\Omega_1)\rightrightarrows(M_1,\pi_1)$ and $(\G_2,\Omega_2)\rightrightarrows(M_2,\pi_2)$ we have the formula
\begin{equation}\label{eq:freecase}
\mu_{\DH}=\vol_1\cdot\vol_2\cdot\mu_{\aff},
\end{equation}
where $\vol_i$ is the $\vol$-function associated to $(\G_i,\Omega_i)\rightrightarrows(M_i,\pi_i)$. Note that as in Remark \ref{rmk:morita} $\mu_{\DH}$ and $\mu_{\aff}$ can be defined using either groupoid and/or action without change. In the free case $\vol_2$ is just the function $\vol_{\red}$ and thus our Main Theorem reduces to (\ref{eq:freecase}).
\end{ex}

\begin{ex}
Consider the special case of $\G$ acting on itself by left translation. This brings us to the situation in Example \ref{ex:freecase}, where the quotient is just $(M,\pi)$ and the integration is $(\G,\Omega)$ (with the same groupoid structure). The Duistermaat-Heckman measure is now the one as defined in \cite[Section 6.3]{pmct2} and equation (\ref{eq:freecase}) becomes
\[ \mu_{\DH}=\vol^2\cdot\mu_{\aff}
\]
which is exactly \cite[Theorem 6.3.1]{pmct2}.
\end{ex}

\begin{ex}
Consider now the more general case of Example \ref{ex:classcase}, namely a locally free Hamiltonian action $G\circlearrowright(\twiddle{X},\twiddle{\omega})$ of a compact, connected Lie group $G$ with proper moment map $\twiddle{\mu}:\twiddle{X}\to\g^*$. As described in Example \ref{ex:classham}, the action lifts to a Hamiltonian action of $(G\ltimes\g^*,-\omega_{\mathrm{can}})$, but this groupoid is not regular unless $G$ is abelian. We can however still consider the regular part $M:=\g^*_{\mathrm{reg}}$: the integration $(\G,\Omega):=(G\ltimes\g_{\mathrm{reg}}^*,-\omega_{\mathrm{can}})$ and its induced action on $(X,\omega):=(\mu^{-1}(M),\twiddle{\omega}|_{\mu^{-1}(M)})$ along $\mu:=\twiddle{\mu}|_X$ inherits all the relevant properties.

Let us now fix a maximal torus $\T\subset G$ and a Weyl chamber $\mathfrak{c}\subset\mathfrak{t}^*$. The coadjoint action gives an isomorphism $(G/\T)\times\mathfrak{c}\cong\g_{\mathrm{reg}}^*$ under which the coadjoint orbit at $\xi\in\mathfrak{c}\subset\g_{\mathrm{reg}}^*$ corresponds to $G/\T\times\{\xi\}$ with the invariant symplectic form $\omega_{\xi}$ which at $e\T$ is given by
\begin{equation}\label{eq:KKSform}
\omega_{\xi}(a,b) = - \xi([a,b]), \,\,\, a,b\in\g/\mathfrak{t}.
\end{equation}
This shows that the Weyl chamber $\mathfrak{c}$ is the leaf space of $\g_{\mathrm{reg}}^*$. This space is smooth, and the integral affine structure is the one dual to the kernel of the exponential map $\exp:\mathfrak{t}\to\T$. In other words, $\mu_{\aff}$ is the standard Lebesgue measure on $\mathfrak{c}$. The function $\vol$ is simply given by the volumes of the coadjoint orbits: the polynomial nature is now immediately clear from Equation (\ref{eq:KKSform}). The function $\vol_{\red}$ is given by the volumes of the original reduced spaces. In view of the discussion in Section \ref{subsec:linvar}, note that the copies $\{g\T\}\times\mathfrak{c}$ provide slices to the coadjoint orbits and that $\T$ is the isotropy group for all points in $\mathfrak{c}$.

As a concrete example, consider $U(k)$ acting on $(\C^{k\times n},\omega_{\mathrm{can}})$ by left multiplication. Identifying $\mathfrak{u}(k)^*\cong\mathfrak{u}(k)$ using the invariant inner product
\[ (A,B)=\mathrm{tr}(A^*B),
\]
this action is Hamiltonian with moment map $\mu:\C^{k\times n}\to\mathfrak{u}(n)$ given by
\[ A\mapsto\frac{i}2AA^*.
\]
A maximal torus $\T\subset U(k)$ is given by the subgroup of diagonal matrices. Its Lie algebra $\mathfrak{t}$ consists of diagonal matrices with purely imaginary entries and we identify it with $\R^k$ in the obvious way. The regular part $\mathfrak{t}_{\mathrm{reg}}$ consists of those matrices with non-repeating entries and we pick the usual fundamental Weyl chamber
\[ \mathfrak{c}=\{p=(p_1,\ldots,p_k)\in\R^k\mid p_1>p_2>\cdots>p_k\}.
\]
For $p\in\mathfrak{c}$, $\mu^{-1}(p)$ is then the ``rescaled'' Stiefel manifold
\[ V_k(\C^n,p)=\{\text{orthogonal } k\text{-frames } (v_1,\ldots, v_k) \text{ in } \C^n \text{ with } \|v_j\|^2 = p_j\}
\]
and the reduced space $\C^{k\times n}\sslash_pU(k)$ is the quotient of $V_k(\C^n,p)$ under the action of $\T=U(1)\times\cdots\times U(1)$ given by
\[ (\lambda_1,\ldots,\lambda_k)\cdot(v_1,\ldots, v_k)=(\lambda_1v_1,\ldots, \lambda_kv_k).
\]
These types of spaces are encountered in estimation theory and optimisation problems on Riemannian quotient manifolds. Our theorem implies that their volumes vary polynomially with $p\in\R^k$.
\end{ex}

\subsection{The proof}\label{subsec:proof}

We now prove the Main Theorem. Let us first address the polynomial nature of $\vol$ and $\vol_{\red}$. For $\vol$ this is well-known (see \cite[Theorem 6.3.1]{pmct2}), and for $\vol_{\red}$ this follows from the discussion in Section \ref{subsec:linvar}: indeed, $\vol_{\red}$ gives precisely the volumes of the reduced spaces of the classical Hamiltonian torus action in the transverse direction and thus it is a polynomial by \cite[Corollary 3.3]{dh}.

It remains to prove equation (\ref{eq:mainthm}). This can be done pointwise, so let us fix $x\in X$ and write $p=\mu(x)$. From the definitions of the affine measure and the Duistermaat-Heckman measure (see Sections \ref{subsec:affmeas} and \ref{subsec:dhmeas} respectively) and equations (\ref{eq:measurecomp}) and (\ref{eq:redspacevol}) it follows that we need to prove the formula
\begin{equation}\label{eq:toprove}
\left|\frac{\omega^{\mathrm{top}}}{\mathrm{top}!}\right|_x=\rho_{\mathrm{Haar}} \otimes \left|\frac{\big(\omega|_{\mu^{-1}(p)}\big)^{\mathrm{top}}}{\mathrm{top}!}\right|_x\otimes(d_x\mu)^*\left(\left|\frac{\omega_{\F_{\pi}}^{\mathrm{top}}}{\mathrm{top}!}\right|_p\otimes|\lambda_1\wedge\cdots\wedge\lambda_q|\right),
\end{equation}
where $\rho_{\mathrm{Haar}}$ is the Haar measure on $\G_p$ (normalised according to $\G_p^0$), viewed as a density on $T_x(\G_p\cdot x)$, and $\{\lambda_1,\ldots,\lambda_q\}$ is a frame of the lattice $\Lambda_p\subset\nu_p^*(\F_{\pi})$.

We define a decomposition $T_xX=V_1\oplus V_2\oplus V_3\oplus V_4$ into linear subspaces consistent with the decomposition in equation (\ref{eq:toprove}) as follows. We write $m=\dim(M)$ and $n=\dim(X)$. Recall also that $q=\dim(\G_p)$. Let $\{\alpha_1,\ldots,\alpha_q\}$ be the basis of $\g_p$ that, through the identifications of Section \ref{subsec:ias}, is equal to $\{\lambda_1,\ldots,\lambda_q\}$, and complete it to a basis $\{\alpha_1,\ldots,\alpha_m\}$ of $A_p$. The infinitesimal action gives us associated vectors $\{\alpha^X_i\}_i$ which form a basis of $T_x(\G\cdot x)$. Next, let $\{v_1,\ldots,v_{n-m-q}\}$ be a basis for a complement of $T_x(\G_p\cdot x)$ in $\ker(d_x\mu)$. Finally, let $\{w_1,\ldots,w_q\}$ be dual to $\{\mu^*\lambda_1,\ldots,\mu^*\lambda_q\}$. Then we define
\begin{align*}
V_1&=\mathrm{span}\{\alpha^X_1,\ldots,\alpha^X_q\} ,\\
V_2&=\mathrm{span}\{v_1,\ldots,v_{n-m-q}\},\\
V_3&=\mathrm{span}\{\alpha^X_{q+1},\ldots,\alpha^X_m\}  ,\\
V_4&=\mathrm{span}\{w_1,\ldots,w_q\} .
\end{align*}

Note that $V_1\oplus V_2=\ker(d_x\mu)$ and $V_1\oplus V_3=\G\cdot T_xX$. Since $\ker(d_x\mu)^{\omega}=\G\cdot T_xX$, it follows that for any $1\leq i \leq q$ the $1$-form $i_{\alpha^X_i}\omega_x$ is automatically zero on $V_1\oplus V_2\oplus V_3$, while on $V_4$ it takes the form
\begin{equation}\label{eq:proof1} \omega_x(\alpha^X_i,w_j)=(\mu^*\lambda_i)(w_j)=\delta_{ij}
\end{equation}
using the moment map condition (\ref{eq:infimoment}). Note also that by definition of the Haar measure we have
\begin{equation}\label{eq:proof2} \rho_{\mathrm{Haar}}(\alpha^X_1,\ldots,\alpha^X_q)=1.
\end{equation}

Next, let us analyse $\omega_x|_{V_3}$. We compute
\begin{align*} \omega_x(\alpha^X_{q+i},\alpha^X_{q+j})&=-\Omega_{1_{\mu(x)}}(\alpha_{q+i},d_x\mu(\alpha^X_{q+j})) \\
&=-\Omega_{1_{\mu(x)}}(\alpha_{q+i},d_{1_{\mu(x)}}\s(\alpha_{q+j})) \\
&=-\Omega_{1_{\mu(x)}}(\alpha_{q+i},\alpha_{q+j}),
\end{align*}
where we use that the $\s$- and $\t$-fibers are $\Omega$-orthogonal. On the other hand, we have
\begin{align*}
(\mu^*\omega_{\F_{\pi}})_x(\alpha^X_{q+i},\alpha^X_{q+j})&=(\omega_{\F_{\pi}})_{\mu(x)}(d_x\mu(\alpha^X_{q+i}),d_x\mu(\alpha^X_{q+j})) \\
&=(\omega_{\F_{\pi}})_{\mu(x)}(d_{1_{\mu(x)}}\s(\alpha_{q+i}),d_{1_{\mu(x)}}\s(\alpha_{q+j})) \\
&=\Omega_{1_{\mu(x)}}(\alpha_{q+i},\alpha_{q+j})
\end{align*}
so that we can conclude that
\begin{equation}\label{eq:proof3} \omega_x(\alpha^X_{q+i},\alpha^X_{q+j})=-(\mu^*\omega_{\F_{\pi}})_x(\alpha^X_{q+i},\alpha^X_{q+j}).
\end{equation}

Combining equations (\ref{eq:proof1}), (\ref{eq:proof2}) and (\ref{eq:proof3}) with the relation $\ker(d_x\mu)^{\omega}=\G\cdot T_xX$ and doing the necessary combinatorics we arrive precisely at equation (\ref{eq:toprove}).

\begin{rmk}
In \cite{Mol2024Stratification} a local model is given for Hamiltonian actions of proper symplectic groupoids, which can be used to give an alternate proof of equation (\ref{eq:toprove}).
\end{rmk}

\bibliographystyle{alpha}

\bibliography{references}

\begin{thebibliography}{BCWZ04}

\bibitem[AB84]{AtiyahBott84Localisation}
M.F. Atiyah and R.~Bott.
\newblock The moment map and equivariant cohomology.
\newblock {\em Topology}, 23(1):1 -- 28, 1984.

\bibitem[Ale97]{Alekseev1997PLGHamAction}
Anton Alekseev.
\newblock On {P}oisson actions of compact {L}ie groups on symplectic manifolds.
\newblock {\em J. Differential Geom.}, 45(2):241--256, 1997.

\bibitem[AMM98]{ALekseevMalkinMeinrenken1998GroupValuedMomentMap}
Anton Alekseev, Anton Malkin, and Eckhard Meinrenken.
\newblock Lie group valued moment maps.
\newblock {\em J. Differential Geom.}, 48(3):445--495, 1998.

\bibitem[Ati82]{Atiyah1982Convexity}
M.~F. Atiyah.
\newblock Convexity and commuting {H}amiltonians.
\newblock {\em Bull. London Math. Soc.}, 14(1):1--15, 1982.

\bibitem[BCWZ04]{BursztynCrainicWeinsteinZhu2004IntDirac}
Henrique Bursztyn, Marius Crainic, Alan Weinstein, and Chenchang Zhu.
\newblock Integration of twisted {D}irac brackets.
\newblock {\em Duke Math. J.}, 123(3):549--607, 2004.

\bibitem[BV82]{BerlineVergne82Localisation}
N.~Berline and M.~Vergne.
\newblock Classes caract{\'e}ristiques {\'e}quivariantes. formule de
  localisation en cohomologie {\'e}quivariante.
\newblock {\em CR Acad. Sci. Paris}, 295(2):539--541, 1982.

\bibitem[CDW87]{CosteDazordWeinstein1987SymmplecticGroupoids}
A.~Coste, P.~Dazord, and A.~Weinstein.
\newblock Groupo\"ides symplectiques.
\newblock In {\em Publications du {D}\'epartement de {M}ath\'ematiques.
  {N}ouvelle {S}\'erie. {A}, {V}ol.\ 2}, volume 87-2 of {\em Publ. D\'ep. Math.
  Nouvelle S\'er. A}, pages i--ii, 1--62. Univ. Claude-Bernard, Lyon, 1987.

\bibitem[CFM19]{pmct2}
M.~Crainic, R.~L. Fernandes, and D.~{Mart{\'i}nez Torres}.
\newblock Regular {P}oisson manifolds of compact types ({PMCT} 2).
\newblock {\em Asterisque}, 413:1--166, January 2019.

\bibitem[CFM21]{LecturesPoisson2021}
M.~Crainic, R.~L. Fernandes, and I.~M\u{a}rcu\c{t}.
\newblock {\em Lectures on {P}oisson geometry}, volume 217 of {\em Graduate
  Studies in Mathematics}.
\newblock American Mathematical Society, Providence, RI, [2021] \copyright
  2021.

\bibitem[CM19]{CrainicMestre2019Measures}
Marius Crainic and Jo\~{a}o~Nuno Mestre.
\newblock Measures on differentiable stacks.
\newblock {\em J. Noncommut. Geom.}, 13(4):1235--1294, 2019.

\bibitem[DH82]{dh}
J.~J. Duistermaat and G.~J. Heckman.
\newblock On the variation in the cohomology of the symplectic form of the
  reduced phase space.
\newblock {\em Inventiones mathematicae}, 69(2):259--268, 1982.

\bibitem[FM24]{FernandesMol2024KahlerMetric}
Rui~Loja Fernandes and Maarten Mol.
\newblock K\"ahler metrics and toric lagrangian fibrations, 2024.
\newblock arXiv 2401.02910.

\bibitem[GS82]{GuilleminSternberg1982Convexity}
V.~Guillemin and S.~Sternberg.
\newblock Convexity properties of the moment mapping.
\newblock {\em Invent. Math.}, 67(3):491--513, 1982.

\bibitem[Kar86]{Karasev1986SymplecticGroupoids}
M.~V. Karas\"ev.
\newblock Analogues of objects of the theory of {L}ie groups for nonlinear
  {P}oisson brackets.
\newblock {\em Izv. Akad. Nauk SSSR Ser. Mat.}, 50(3):508--538, 638, 1986.

\bibitem[Lu91]{Lu1989PLGMomentMap}
Jiang-Hua Lu.
\newblock Momentum mappings and reduction of {P}oisson actions.
\newblock In {\em Symplectic geometry, groupoids, and integrable systems
  ({B}erkeley, {CA}, 1989)}, volume~20 of {\em Math. Sci. Res. Inst. Publ.},
  pages 209--226. Springer, New York, 1991.

\bibitem[McD88]{McDuff1988QuasiHam}
Dusa McDuff.
\newblock The moment map for circle actions on symplectic manifolds.
\newblock {\em J. Geom. Phys.}, 5(2):149--160, 1988.

\bibitem[Mey73]{Meyer1973SympRed}
K.~R. Meyer.
\newblock Symmetries and integrals in mechanics.
\newblock In {\em Dynamical systems ({P}roc. {S}ympos., {U}niv. {B}ahia,
  {S}alvador, 1971)}, pages 259--272. Academic Press, New York-London, 1973.

\bibitem[Mol23]{Mol2023MultFree}
Maarten Mol.
\newblock On the classification of multiplicity-free hamiltonian actions by
  regular proper symplectic groupoids, 2023.
\newblock arXiv 2401.00570.

\bibitem[Mol24]{Mol2024Stratification}
Maarten Mol.
\newblock Stratification of the transverse momentum map.
\newblock {\em Compositio Mathematica}, 160(3):518–585, 2024.

\bibitem[MW74]{MarsdenWeinstein1974SympRed}
J.~Marsden and A.~Weinstein.
\newblock Reduction of symplectic manifolds with symmetry.
\newblock {\em Reports on Mathematical Physics}, 5(1):121--130, 1974.

\bibitem[MW88]{MikamiWeinstein1988HamiltonianGroupoidActions}
Kentaro Mikami and Alan Weinstein.
\newblock Moments and reduction for symplectic groupoids.
\newblock {\em Publ. Res. Inst. Math. Sci.}, 24(1):121--140, 1988.

\bibitem[SL91]{SjamaarLerman1991SingRed}
R.~Sjamaar and E.~Lerman.
\newblock Stratified symplectic spaces and reduction.
\newblock {\em Annals of Mathematics}, 134(2):375--422, 1991.

\bibitem[Wei87]{Weinstein1987SymplecticGroupoids}
Alan Weinstein.
\newblock Symplectic groupoids and {P}oisson manifolds.
\newblock {\em Bull. Amer. Math. Soc. (N.S.)}, 16(1):101--104, 1987.

\bibitem[WX92]{WeinsteinXu1992QYBE}
Alan Weinstein and Ping Xu.
\newblock Classical solutions of the quantum {Y}ang-{B}axter equation.
\newblock {\em Comm. Math. Phys.}, 148(2):309--343, 1992.

\bibitem[Xu04]{Xu2004MomentumMorita}
Ping Xu.
\newblock Momentum maps and {M}orita equivalence.
\newblock {\em J. Differential Geom.}, 67(2):289--333, 2004.

\end{thebibliography}

\end{document}